# Euler constant as a renormalized value of Riemann zeta function at its pole. Rationals related to Dirichlet *L*-functions

Andrei Vieru

**Abstract**

We believe that Euler's constant is not just the "renormalized" value of the Riemann zeta function in 1. In a sense that we shall clarify it is in fact the ***normal and natural*** value of zeta of 1. In this paper we first propose a limit definition of a function whose values coincide everywhere with those of the Riemann zeta function, save in 1, where our limit definition yields the Euler constant. Since in the literature one can find more than one way to regularize the value of the zeta function at $s=1$, we give asymptotic expansions where, by dint of some extended analogies, Euler's constant appears to be the true "renormalized" value. As a striking example of such analogies, we propose an expansion of the logarithm function based on Euler's constant and on all values of the zeta function at odd positive integers, in which all these presumably irrational numbers are accompanied by Harmonic numbers of corresponding orders. The other aim of this paper is to show how sequences of rationals, often the same, arise in computations related to Dirichlet *L*-functions. Here, a connection with the Liouville lambda function appears to have been found. Thus we raise the question about the possible usefulness of an extension of the Liouville lambda function to rationals.

## 1. Motivation

Euler's constant viewed as a renormalized value of the Riemann zeta function at $s = 1$ is still an 'active area of current research' since Stephen Hawking has published in 1977 *Zeta function regularization of path integrals in curved spacetime, (Comm. Math. Phys. 55)*.

We give an infinite family of expansion formulae of $\ln\Gamma$ none of which are around $z = 1$, and whose summarizing into one single formula is nevertheless much more *illustrative* for our issue than just the convention of taking as a renormalized value the constant term in the Laurent expansion of the Riemann zeta function near $s = 1$. This convention (as pointed out by Jeffrey C. Lagarias) 'also matches the constant term in the Taylor expansion of the digamma function $\psi(z)$ around $z = 1$' (see [11]).

For the sections 4-7, the real motivation is the following one: since not only we know the distribution of the first billions of primes, but we know them all, so to say, by their names, it seemed interesting to try to grasp something about primes starting the other way around, 'way back from infinity'. Since, due to the Euler product formula, Riemann zeta function is about primes not just in the critical strip, but virtually everywhere, we tried to compute asymptotic expansions of functions of the form $1/(L(s) - 1)$ in the neighborhoods of infinity, where $L$ is either the Riemann zeta function itself or some other Dirichlet $L$-function.

To describe the irregular features of these expansions, Legendre, Jacobi and Kronecker symbols seemed of no avail. The sequences of irregular signs and other irregular features of the rationals that arise in all expansions of this type can be partly described in terms of Liouville's lambda function (extended in this context to the class of rationals whose numerators are integer powers of primes). We formulate a not very obvious connection between these expansions and the Liouville lambda function, which is known to be directly related both to primes and to Riemann zeta function. Higher accuracy computation capabilities are needed to get more terms in these asymptotic expansion. These expansions, where $s$ and the fast growing functions $1/(L(s) - 1)$ both run to infinity, are of interest: they are made of either surprisingly small primes or larger composite integers with still surprisingly small prime factors, very often just integer powers of small primes. One can say that the whole study is about the way small numbers — particularly small primes or products of small primes — are mirrored in an 'infinity' related to the Riemann zeta function and to other Dirichlet $L$-functions. The tenuous link between the two themes treated in this paper lies in the possibility to express renormalized values of $L$-functions at $s = 1$ in terms of iterated $L$-functions and of powers of rationals.

## 2.0. Why the value of zeta of 1 may be viewed as a finite real number

The function that has been named later the 'Riemann zeta function' was first defined by Euler in the real domain by means of a series that converges for $s > 1$:

$$\zeta(s) = \sum_{n=1}^{\infty} \frac{1}{n^s} \qquad (1°)$$

It is indeed possible to define for $s \geq 1$ the following family of functions:

$$\zeta_r(s) = \lim_{k \to \infty} \left\{ \sum_{n=1}^{k} \frac{1}{n^s} - \sum_{m=k+1}^{\lfloor rk^2 \rfloor} \frac{1}{m^s} \right\} \qquad (2°)$$

where $r$ is a strictly positive real number, and where $\lfloor x \rfloor$ is the floor function.

Clearly, for any $r > 0$, and for any real $s > 1$, $\quad \zeta_r(s) = \sum_{n=1}^{\infty} \frac{1}{n^s} = \zeta(s) \quad$ **(3°)**

(since, however small a given $r$ may be, for sufficiently large $k$, we'll have $\lfloor rk^2 \rfloor > k+1$, and since the subtracted sum tends to 0 when $k \to \infty$)

It is easy to see that $\quad \zeta_r(1) = \gamma - \ln r \quad$ **(4°)**

This value becomes obvious if one thinks, geometrically, that for sufficiently large $k$ the approximate equality

$$\sum_{m=k+1}^{\lfloor rk^2 \rfloor} \frac{1}{m} \approx \int_{k}^{\lfloor rk^2 \rfloor} \frac{dx}{x} \approx \ln r + \ln k^2 - \ln k = \ln r - \ln k \quad \textbf{(5°)}$$

tends to a genuine equality when $k \to \infty$.

Substituting $\ln r - \ln k$ into **(2°)** (in which we set $s = 1$) and knowing that Euler's constant is defined as

$$\gamma = \lim_{k \to \infty} \left\{ \sum_{n=1}^{k} \frac{1}{n} - \ln k \right\} \quad \textbf{(6°)}$$

we get **(3°)**.

In particular[1], if we set $r = 1$ we have $\quad \zeta_1(1) = \gamma \quad$ **(7°)**

## 2.1 Why zeta of 1 cannot be anything else but Euler's constant

We've seen that a far-fetched but not incorrect way to define $\varsigma(s)$ is the following one:

$$\zeta(s) = \lim_{k \to \infty} \left\{ \sum_{n=1}^{k} \frac{1}{n^s} - \sum_{m=k+1}^{k^2} \frac{1}{m^s} \right\} \quad \textbf{(8°)}$$

Ok, that's an arbitrary far-fetched choice, but then, according to a formula of Macys, and to what has been said in the previous chapter

$$\zeta(1) = \lim_{k \to \infty} \left\{ \sum_{j=1}^{k} \frac{1}{j} - \sum_{l=k+1}^{k^2} \frac{1}{l} \right\} = \gamma \quad \textbf{(9°)}$$

---

[1] To get $\varsigma(1) = 0$ as some others use to "renormalize" the value of zeta at $s = 1$, we have to set $r = e^{\gamma}$ which may be "interesting" but less "natural".

(Here again, and clearer than ever, the integral representation of ln*k* as $\int_{k}^{k^2} \frac{1}{x} dx$ might be the origin of a simple geometric proof of Macys' formula.)

But that's not enough to establish that zeta of 1 equals Euler's constant. We need to see why, as we did in **(7°)**, we should set *r* = 1 except that it is the natural and the simplest one. Actually, there are several well-known Taylor series where one may see that Euler's constant 'behaves' as if it where zeta of 1. The more complicated will be the found formulae[1], the more convincing will appear the thesis of the Euler's constant as being the "true" regularized value of the Riemann zeta function in the point which appears as an "unremovable singularity".

## 2.2 Expanding lnΓ near the negative poles of the Gamma function

The digamma function Taylor expansion that Professor Jeffrey C. Lagarias refers to is

$$\psi(z+1) = -\gamma + \sum_{k=1}^{\infty} (-1)^{k+1} \zeta(k+1) z^k \qquad \text{(for } |z| < 1\text{)} \qquad \textbf{(10°)}$$

Here the only accompanying element that fits the analogy is the (negative) sign.

Near 0, lnΓ may be written as

$$\ln\Gamma(x) = -\ln x - \gamma x + \sum_{k=2}^{\infty} \frac{\zeta(k)}{k}(-x)^k \qquad \text{(for } x < 1\text{)} \qquad \textbf{(11°)}$$

If $\zeta(1) = \gamma$ then in the RHS it would be possible to drag the second summand under the infinite summation, thus beginning it from *k* = 1. Here the analogy slightly broader: the denominators fit as well as the signs of the summand.

To make more convincing and completely unambiguous this analogy between γ and the values of the Riemann zeta function at integer arguments greater than 1, we shall now expand lnΓ near the negative poles of the Gamma function.

---

[1] with the more encompassing parts of them fitting some obvious global analogy (necessarily including the analogy between Euler's constant and the values of Riemann zeta function for greater integers)

Bellow we shall provide an expansion in which this analogy becomes much more obvious, since it is established not just between Euler's constant and the values of the zeta function at greater integer arguments, but between larger parts of the considered formula. In our formula, as the index $k$ runs through all integers, the summands $\varsigma(k)$ appear escorted by harmonic numbers of order $k$ while $\gamma$ appears in company of ordinary harmonic numbers. One will have to notice that the analogy we are speaking about is strengthened by the fact that it holds for all non zero negative integers where the Gamma function has singularities.

In the real domain, one has for all positive integers $n$:

$$\ln\Gamma\bigl(-n+(-1)^n x\bigr) = -\ln x - \sum_{j=1}^{n} \ln j +$$

$$+(-1)^{n+1}(\gamma - H_n)x + \sum_{k=2}^{\infty} (-1)^{(n+1)k} \frac{\varsigma(k)+(-1)^k H_{n,k}}{k} x^k \qquad (12°)$$

(where $H_{n,k}$ are harmonic numbers of rank $n$ and of order $k$)

Here again, under the convention $\varsigma(1) = \gamma$, in the RHS of **(12°)** one can 'drag' the third summand under the infinite summation, which, doing so, would start with $k = 1$

$$\ln\Gamma\bigl(-n+(-1)^n x\bigr) = -\ln x - \sum_{j=1}^{n} \ln j + \sum_{k=1}^{\infty} (-1)^{k(n+1)} \frac{\varsigma(k)+(-1)^k H_{n,k}}{k} x^k$$

$$(13°)$$

Note also that ***in the numerator***, as $k$ runs through the natural integers, the signs of the second summand alternate for all $n$, while the signs of the terms of the ***whole infinite sum*** alternate only for even $n$. Furthermore, the analogy between Euler's constant viewed as $\varsigma(1)$ and the values of $\varsigma$ at greater integer arguments is reinforced by the presence of the exponent of $-1$, which depends on both $k$ and $n$. Last but not least, as in **(11°)** the $k$ in the denominator completes the domain of the analogy.

(We consider $\ln\Gamma$ only in the real domain, avoiding possible but unnecessary discussions in this context about branches in the complex domain inherited from the logarithmic function.)

## 2.3. The logarithm function: an asymptotic expansion based on Euler's constant and on the values of the Riemann zeta function at odd positive integers

In the previous section, we have pointed out how Euler's constant may convincingly show up as the renormalized value of the Riemann zeta function at its pole in expansions of some special functions. Now we'll show how Euler's constant may reveal itself as the renormalized value of the zeta function in infinite sums representing mere *numbers*. The expansion

$$\ln\frac{3}{2} = -(\gamma - 1) - \frac{1}{2^2}\frac{\zeta(3)-1}{3} - \frac{1}{2^4}\frac{\zeta(5)-1}{5} - \ldots = -\sum_{k=1}^{\infty}\frac{1}{2^{2(k-1)}}\frac{\zeta(2k-1)-1}{2k-1}$$

(we use the convention $\varsigma(1) = \gamma$) was known to Euler and later to Stieltjes; they used this formula to compute $\gamma$ up to the 12th, and respectively 33rd, decimal.

On the other hand, **(13°)** suggests that $\ln\Gamma(-3/2)$ might be written as a power expansion of $x = 1/2$ in two different ways: once considering $-3/2$ as belonging to the neighborhoods of $-1$, and, as well, in a different form, considering it as belonging to the vicinities of $-2$. This reduces to set in **(13°)** $n=1$, respectively $n=2$

It turns out that both expansions are absolutely correct. Subtracting one expansion from the other and using some elementary algebra, one gets

$$\ln 2 = -\frac{2\gamma - 2\frac{1}{2}}{1}\frac{1}{2} + \frac{1}{8}\frac{1}{2^2} - \frac{2\zeta(3) - 2\frac{1}{8}}{3}\frac{1}{2^3} + \frac{1}{64}\frac{1}{2^4} - \frac{2\zeta(5) - 2\frac{1}{32}}{5}\frac{1}{2^5} + \ldots$$

This elementary technique easily generalizes to

$$\ln(2n) = -\sum_{k=1}^{\infty}\frac{2\zeta(2k-1) - H_{2n-1}^{(2k-1)} - H_{2n}^{(2k-1)}}{2k-1}\frac{1}{2^{2k-1}} + \sum_{j=1}^{\infty}\frac{1}{2j(2n)^{2j}}\frac{1}{2^{2j}}$$

where $H_j^{(m)}$ – usually written as $H_{j,\,m}$ – are the generalized harmonic numbers of order $m$ and argument $j$, and where $\zeta(2k-1)$ is taken to be $\gamma$ when $k = 1$.

This is another nice appearance of Euler's constant as a renormalized value of the Riemann zeta function at its pole.

One may believe things might be trickier with the intervals of the real axis where the Gamma function takes negative values. Actually, for our purpose, one doesn't need to know how to compute logarithms of negative real numbers. Using the same technique — i.e. writing $\ln\Gamma(-2n-1/2)$ in both ways, as if **(13°)** would apply to neighborhoods in either direction — one immediately finds that the first term in **(13°)**, namely $\ln(-1/2)$, is always canceled. So, under the same conventions, we have

$$\ln n = \sum_{k=1}^{\infty} \left( \frac{1}{2kn^{2k}} \frac{1}{2^{2k}} - \frac{2\zeta(2k-1) - H_{n-1}^{(2k-1)} - H_n^{(2k-1)}}{2k-1} \frac{1}{2^{2k-1}} \right) \quad \textbf{(14°)}$$

However, if we were to send this formula to some knowledgeable mathematical website, we would probably drop the convention $\varsigma(1) = \gamma$, and write **(14°)** in a more traditional manner, namely as

$$\ln n = \frac{1}{8n^2} - \frac{2\gamma - H_{n-1} - H_n}{2}$$

$$+ \sum_{k=2}^{\infty} \left( \frac{1}{2kn^{2k}} \frac{1}{2^{2k}} - \frac{2\zeta(2k-1) - H_{n-1}^{(2k-1)} - H_n^{(2k-1)}}{2k-1} \frac{1}{2^{2k-1}} \right) \quad \textbf{(14a°)}$$

This formula holds as well for non integer arguments greater than 1 using (generalized) Harmonic numbers with fractional argument[1]. These are introduced through:

$$H_{q/p,n} = \zeta(n) - p^n \sum_{k=1}^{\infty} \frac{1}{(q+pk)^n}$$

There exists a classical expansion:

$$\ln n = -\gamma + H_n - \frac{1}{2n} + \sum_{k=1}^{\infty} \frac{B_{2k}}{2kn^{2k}} \quad \textbf{(15°)}$$

which also works for fractional *n*. It might be compared with ours, inasmuch Euler's constant, Harmonic numbers and the form $2kn^{2k}$ appear in both of them. Zeta values for odd arguments appear under the infinite sum only in our expansion, while in **(15°)** we see Bernoulli numbers (which, instead, show up in the closed-form expression of the Riemann zeta function at integer even arguments).

As we already mentioned it, our purpose was not to find new expansions for the sake of the logarithm function, but for the sake of Euler's constant ***as a regularized value of Riemann zeta function at its pole.***

In our mind, the discontinuity of the Riemann zeta function in 1 is, in some sense, analogous to the discontinuity in the family of antiderivatives $\int dx/x^\alpha$ when *α*=1.

---

[1] Note that Euler's constant may be expanded without any use of logarithms (see [12])

## 2.4. Euler's constant and ζ(2) expressed in terms of Log Γ

In the general case, for non-negative *n*, we have:

$$\lim_{x \to 0} \left\{ \frac{\ln(n!x) + \ln\Gamma(-n + (-1)^n x)}{x} \right\} = (-1)^{n+1}(\gamma - H_n) \quad \textbf{(A°)}$$

Replacing $(-1)^n$ and ln –Γ by ln|Γ(*x*)| we'll have

$$\lim_{x \to 0} \left\{ \frac{\ln(n!x) + \ln|\Gamma(-n \pm x)|}{x} \right\} = \pm(\gamma - H_n) \quad \textbf{(B°)}$$

The RHS of this limit identity might be written as

$$\pm\left( ''\zeta(1)'' - \sum_{k=1}^{n} \frac{1}{k} \right) \quad \textbf{(C°)}$$

(Are quotation marks still needed, because thinking of Euler's constant as "ζ(1)" is, in some sense, still a metaphor? As generally known, and as already said, according to its original definition the Riemann zeta function has a singularity in 1, which is not removable, and which has residue 1. The philosophical and semiotical problem of metaphors in mathematics, and more generally in science, will not be discussed here.) One can notice that although the formulae

$$\lim_{x \to 0} \left\{ \frac{\ln(n!x) + \ln|\Gamma(-n + x)|}{x} \right\} = -(\gamma - H_n) \quad \textbf{(D°)}$$

and

$$\lim_{x \to 0} \left\{ \frac{\ln(n!x) + \ln|\Gamma(-n - x)|}{x} \right\} = +(\gamma - H_n) \quad \textbf{(E°)}$$

state convergence to numbers which are identical in absolute value, an arbitrarily chosen sequence of *x* that converges to 0, does not yield identical sequences (in absolute value) for the expressions in the LHS of **(D°)** and of **(E°)**.

One can effectively 'compare' these sequences adding the expressions that stand under the limit in **(D°)** and in **(E°)** and then dividing the result by *x*.

Doing so, for *n* = 0, one can find

$$\lim_{x \to 0} \left\{ \frac{2\ln x + \ln|\Gamma(-x)| + \ln|\Gamma(x)|}{x^2} \right\} = \zeta(2) = \frac{\pi^2}{6} \qquad \textbf{(F°)}$$

and, for *n* = 1,

$$\lim_{x \to 0} \left\{ \frac{2\ln x + \ln|\Gamma(-1-x)| + \ln|\Gamma(-1+x)|}{x^2} \right\} = \zeta(2) + 1 \qquad \textbf{(G°)}$$

In the general case, one has:

$$\lim_{x \to 0} \left\{ \frac{2\ln(n!x) + \ln|\Gamma(-n-x)| + \ln|\Gamma(-n+x)|}{x^2} \right\} = \zeta(2) + \sum_{k=1}^{n} \frac{1}{k^2} \qquad \textbf{(H°)}$$

and, therefore,

$$\lim_{\substack{x \to 0 \\ n \to \infty}} \left\{ \frac{2\ln(n!x) + \ln|\Gamma(-n-x)| + \ln|\Gamma(-n+x)|}{x^2} \right\} = 2\zeta(2) \qquad \textbf{(I°)}$$

"Comparing" the process of convergence at two consecutive singularities of the gamma function we get:

$$\lim_{x \to 0} \left\{ \frac{1}{x^2} \left[ 1 - \frac{2\ln(n!x) + \ln|\Gamma(-n-x)| + \ln|\Gamma(-n+x)| - x^2 \sum_{k=1}^{n} \frac{1}{k^2}}{2\ln[(n+1)!x] + \ln|\Gamma(-n-1-x)| + \ln|\Gamma(-n-1+x)| - x^2 \sum_{k=1}^{n+1} \frac{1}{k^2}} \right] \right\}$$

$$= \frac{3}{(n+1)^4 \pi^2} = \frac{1}{2(n+1)^4 \zeta(2)}$$

The following limit formula holds as well[1]:

$$\lim_{x \to 0}\left\{\frac{1}{x^2}\ln\left(\frac{2\ln(n!x)+\ln|\Gamma(-n-x)|+\ln|\Gamma(-n+x)|-x^2\sum_{k=1}^{n}\frac{1}{k^2}}{2\ln[(n+1)!x]+\ln|\Gamma(-n-1-x)|+\ln|\Gamma(-n-1+x)|-x^2\sum_{k=1}^{n+1}\frac{1}{k^2}}\right)\right\}$$

$$=-\frac{1}{2(n+1)^4\zeta(2)}$$

## 3.1 Integer powers of ln 2 expressed in terms of ln$\zeta$ and ln$\eta$

$$\lim_{x \to \infty}\left\{\frac{\ln\zeta(\ln\zeta(x)+1)}{x}\right\}=\lim_{x \to \infty}\left\{\frac{\ln\zeta(\zeta(x))}{x}\right\}=\ln 2$$

Let $f_1(x)=\ln\zeta(\ln\zeta(x)+1)$ and $f_2(x)=\ln(\ln\zeta(x))$

then, iterating, we'll have, trivially, for any $n$

$$\lim_{x \to \infty}\left\{\frac{f_1^n(x)}{x}\right\}=\lim_{x \to \infty}\left\{\frac{f_2^n(x)}{x}\right\}=(\ln 2)^n$$

and, somewhat less trivially,

$$\lim_{x \to \infty}\{2^x[\ln\zeta(\ln\zeta(x)+1)-\ln\zeta(\zeta(x))]\}=\frac{1}{2}$$

The Dirichlet eta function is defined as

$$\eta(s)=\sum_{n=1}^{\infty}\frac{(-1)^{n-1}}{n^s}=\frac{1}{1^s}-\frac{1}{2^s}+\frac{1}{3^s}-\frac{1}{4^s}+\ldots$$

a series that converges for any complex number with real part greater than 0.

---

[1] In this limit formula $n!$ and $(n+1)!$ appear as absolute values of the residues in the singularity $-n$ and $-(n+1)$ of the gamma function (see also [3]).

Define $\quad g(x) = |\ln|\ln\eta(x)||$

Then for any $n \geq 0$ $\quad \lim\limits_{x\to\infty}\left\{\dfrac{g^n(x)}{x}\right\} = (\ln 2)^n$

Define $\quad h(x) = \ln|\zeta(\eta(x))|$

then, again (more or less trivially), for any $n \geq 0$ $\quad \lim\limits_{x\leftarrow\infty}\left\{\dfrac{h^n(x)}{x}\right\} = (\ln 2)^n$

## 3.2. The iterations of a pair of functions "preserving asymptotic proportions at infinity"

The nice thing is that for any $N \geq 1$ we have the following general formula, which holds for all iterates of $g$ and $h$ (of strictly positive order):

$$\lim\limits_{x\to\infty}\left\{2^{x(\ln 2)^{n-1}}[g^n(x) - h^n(x)]\right\} = \gamma - \dfrac{1}{2} \qquad (\odot)$$

$$\lim\limits_{x\to\infty}\left\{2^{x-1}[f(x) - h(x)]\right\} = \gamma \qquad (\ast)$$

but no formula is known to us concerning higher order iterates of $f$ and $h$, let alone concerning relations that involve higher order iterates of $f$ and $g$.

Another nice formula involving iterates of $g$ and of $h$ is the following one:

Let $\quad \omega(x, n) = 2^{x(\ln 2)^{n-1}}[g^n(x) - h^n(x)] \qquad (\diamondsuit)$

then again, for any positive integer $n$, we'll have

$$\lim\limits_{x\to\infty}\left\{\dfrac{\omega(x, n+1) - \omega(x, n)}{\ln\omega(x, n+1) - \ln\omega(x, n)}\right\} = \gamma - \dfrac{1}{2} \qquad (\ast)$$

## 4. SEQUENCES OF RATIONALS ARISING IN CONNECTION WITH DIRICHLET *L*-SERIES

Trying to better understand (⊙), (✳), (❖) and (✺), I have begun expanding $1/(\zeta(x) - 1)$ for large values.

Having found that, for sufficiently remote *x* values,

$$\frac{1}{\zeta(x) - 1} = 2^x - \left(\frac{4}{3}\right)^x - 1 + \left(\frac{8}{9}\right)^x - \left(\frac{4}{5}\right)^x + \left(\frac{2}{3}\right)^x - \left(\frac{16}{27}\right)^x - \left(\frac{4}{7}\right)^x +$$

$$+ 2\left(\frac{8}{15}\right)^x - 2\left(\frac{4}{9}\right)^x + \left(\frac{2}{5}\right)^x + \left(\frac{32}{81}\right)^x + 2\left(\frac{8}{21}\right)^x - \left(\frac{4}{11}\right)^x - 3\left(\frac{16}{45}\right)^x + \left(\frac{8}{25}\right)^x -$$

$$-\left(\frac{4}{13}\right)^x + 3\left(\frac{8}{27}\right)^x + \left(\frac{2}{7}\right)^x - 3\left(\frac{4}{15}\right)^x - \left(\frac{64}{243}\right)^x - \ldots \quad (\text{♚})$$

I filled the first five or six denominators in (♚) as an entry search into OEIS, and learned that this expansion was already known to Benoît Cloitre (A112932), to whom I express my admiration for his entire work. (I just added several more terms to his expansion; for even more terms, see below, page 20 of this paper.)

The derived Cloitre's formula

$$\lim_{x \to \infty} \left\{ \zeta(\zeta(x)) - 2^x + \left(\frac{4}{3}\right)^x + 1 \right\} = \gamma \quad (\text{♛})$$

is surprising for who doesn't know (♚). In fact there are uncountable analogues of this limit formula (involving or not Euler's constant).

Our first observation is that several functions have an expansion with, up to the constant term, the same beginning:

$$\frac{1}{\ln \zeta(s)} = 2^s - \left(\frac{4}{3}\right)^s - \frac{1}{2} + \left(\frac{8}{9}\right)^s - \left(\frac{4}{5}\right)^s + \left(\frac{2}{3}\right)^s - \left(\frac{16}{27}\right)^s - \left(\frac{4}{7}\right)^s + 2\left(\frac{8}{15}\right)^s - \ldots$$

$$\frac{1}{\ln(\ln \zeta(x) + 1)} = 2^x - \left(\frac{4}{3}\right)^x - 0 + \left(\frac{8}{9}\right)^x - \left(\frac{4}{5}\right)^x + \left(\frac{2}{3}\right)^x - \left(\frac{16}{27}\right)^x - \left(\frac{4}{7}\right)^x + 2\left(\frac{8}{15}\right)^x - \ldots$$

(with a gain of 1/2 in the constant term for each substitution of *x* by $\ln\zeta(x)+1$)

The next example of substitution (which we restrain from writing down) yields the following analogue of (♛):

$$\lim_{x \to \infty} \left\{ \zeta[\ln[\ln(\ln \zeta(x) + 1) + 1] + 1] - 2^x + \left(\frac{4}{3}\right)^x - \frac{1}{2} \right\} = \gamma$$

The natural question about these expansions is whether other reciprocals of modified Dirichlet *L*-functions permit similar asymptotic expressions in terms of sums of powers of rationals. Here are but a few examples for the most popular Dirichlet eta[1], beta end lambda *L*-functions:

$$\frac{1}{1-\eta(x)} = 2^x + \left(\frac{4}{3}\right)^x - 1 + \left(\frac{8}{9}\right)^x + \left(\frac{4}{5}\right)^x - 3\left(\frac{2}{3}\right)^x + \left(\frac{16}{27}\right)^x +$$

$$+ \left(\frac{4}{7}\right)^x + 2\left(\frac{8}{15}\right)^x - 4\left(\frac{4}{9}\right)^x - 3\left(\frac{2}{5}\right)^x + \ldots \quad (\text{♚})$$

As one can see, in this peculiar case, the rationals are the same, but the coefficients (and in particular the signs) are not. (Further terms in the expansion (♚) of $\eta$ are not similar to those of ς.) The presence of powers of 2 in the numerators could not escape our attention, nor did the presence of powers of 3 in the two following examples:

Let $\beta$ be the Dirichlet beta function, defined as

$$L(s, \chi_{-2}) = \beta(s) = \sum_{n=0}^{\infty} \frac{(-1)^n}{(2n+1)^s}$$

Then, for large *s*

$$\frac{1}{1-\beta(s)} = 3^s + \left(\frac{9}{5}\right)^s - \left(\frac{9}{7}\right)^s + \left(\frac{27}{25}\right)^s + 1 - \left(\frac{9}{11}\right)^s - 2\left(\frac{27}{35}\right)^s +$$

$$+ \left(\frac{9}{13}\right)^s + \left(\frac{81}{125}\right)^s + \left(\frac{3}{5}\right)^s + \left(\frac{27}{49}\right)^s + \left(\frac{9}{17}\right)^s - 2\left(\frac{27}{55}\right)^s - \left(\frac{9}{19}\right)^s + \ldots \quad (\text{♛})$$

Let $\quad L(s, \chi_2) = \lambda(s) \equiv \sum_{n=0}^{\infty} \frac{1}{(2n+1)^s}$

Then, for large *s*

$$\frac{1}{\lambda(s)-1} = 3^s - \left(\frac{9}{5}\right)^s - \left(\frac{9}{7}\right)^s + \left(\frac{27}{25}\right)^s - 1 - \left(\frac{9}{11}\right)^s + 2\left(\frac{27}{35}\right)^s -$$

$$- \left(\frac{9}{13}\right)^s - \left(\frac{81}{125}\right)^s + \left(\frac{3}{5}\right)^s + \left(\frac{27}{49}\right)^s - \left(\frac{9}{17}\right)^s + 2\left(\frac{27}{55}\right)^s - \left(\frac{9}{19}\right)^s - 3\left(\frac{81}{175}\right)^s$$

$$+ \left(\frac{3}{7}\right)^s + 2\left(\frac{27}{65}\right)^s - \left(\frac{9}{23}\right)^s + \left(\frac{243}{625}\right)^s - 2\left(\frac{9}{25}\right)^s + 2\left(\frac{27}{77}\right)^s - 2\left(\frac{81}{245}\right) - \ldots \quad (\text{♖})$$

---

[1] here is only a preview of the very particular Dirichlet $\eta$-function, whose asymptotic expansion will be studied further bellow.

The last expansion yields, the two following limit formulae (and many others):

$$-\lim_{x\to\infty}\left\{\Gamma(\lambda(x)-1)-3^x+\left(\frac{9}{5}\right)^x+\left(\frac{9}{7}\right)^x-\left(\frac{27}{25}\right)^x+1\right\}=\gamma$$

Since $\lim_{x\to 0}\left\{\lambda(1+x)-\frac{1}{2x}\right\}=\frac{1}{2}(\gamma+\ln 2)$

we'll have, using just the first few terms of (♖):

$$\lim_{x\to\infty}\left\{\lambda(\lambda(x))-\frac{1}{2}\left[3^x-\left(\frac{9}{5}\right)^x-\left(\frac{9}{7}\right)^x+\left(\frac{27}{25}\right)^x-1\right]\right\}=\frac{1}{2}(\gamma+\ln 2)$$

The RHS of this formula 'happens' to be one of the Euler-Lehmer constants, namely γ(1,2) (see Lehmer's paper, [10]). One can add more terms from (♖) into the brackets in the LHS of this limit formula and that will speed up the convergence.

Note that $\frac{1}{2}(\gamma+\ln 2)$ is nothing else but the renormalized value of $\lambda$ at $s=1$.

Now the real question is: given a descending (by absolute value) infinite sequence of rationals $(q)_{n\in N}$, given a sequence of integer coefficients $(a_n)$ and given the corresponding infinite series

$$V(x)=a_1 q_1^x + a_2 q_2^x + \ldots$$

*to what class or subclass of functions belongs $1/(V(x)+C)$ ?* (where C is a constant, possibly 0)

The class of these functions may not necessarily coincide with the Dirichlet *L*-functions. Surprisingly, some irregular series give birth to similar sequences of $(q)_{n\in N}$ as the Riemann zeta function does, although not with the same coefficients $(a_n)$.

Here are two examples: $\vartheta(s)=\sum_p p^{-s}$

(the prime Zeta function where the sum is taken over all primes.) We have:

$$\frac{1}{\vartheta(s)}=2^s-\left(\frac{4}{3}\right)^s+\left(\frac{8}{9}\right)^s-\left(\frac{4}{5}\right)^s-\left(\frac{16}{27}\right)^s-\left(\frac{4}{7}\right)^s+2\left(\frac{8}{15}\right)^s+\left(\frac{32}{81}\right)^s$$

$$+2\left(\frac{8}{21}\right)^s-\left(\frac{4}{11}\right)^s-3\left(\frac{16}{45}\right)^s+\left(\frac{8}{25}\right)^s-\left(\frac{4}{13}\right)^s-\left(\frac{64}{243}\right)^s-3\left(\frac{16}{63}\right)^s+2\left(\frac{8}{33}\right)^s+\ldots$$

$$+4\left(\frac{32}{135}\right)^s-\left(\frac{4}{17}\right)^s+2\left(\frac{8}{35}\right)^s-3\left(\frac{16}{75}\right)^s-\left(\frac{4}{19}\right)^s+2\left(\frac{8}{39}\right)^s+\left(\frac{8}{45}\right)^s\ldots \qquad (\knight)$$

Comparing (♞) with (♚), one finds among the initial terms of these asymptotic expansions exactly the same summands just with the 3rd, 6th, 10th, 11th, 18th, 19th, 20th, 27th, 31st, 32nd, 33rd… terms skipped, namely −1, (2/3)^x, −2(4/9)^x, (2/5)^x, 3(8/27)^x, (2/7)^x, −3(4/15)^x, etc…

The nice thing is that (♞) and (♚) seem to obey to a similar law related to the Liuoville lambda function. This issue will be discussed later.

Here are two more example of irregular series:

let $\daleth(s) = \sum_{n \in A} n^{-s} - \sum_{m \in B} m^{-s}$

where $A$ is the set of numbers which correspond to the ranks of zeros in the Thue-Morse sequence, whereas $B$ corresponds to the set of the ranks of 1 in the same sequence.

Then, for large $x$

$$\frac{1}{1 - \daleth(x)} = 2^x - \left(\frac{4}{3}\right)^x + 1 + \left(\frac{8}{9}\right)^x - \left(\frac{4}{5}\right)^x - \left(\frac{2}{3}\right)^x - \left(\frac{16}{27}\right)^x + \left(\frac{4}{7}\right)^x +$$

$$+ 2\left(\frac{8}{15}\right)^x - \left(\frac{2}{5}\right)^x + \left(\frac{32}{81}\right)^x - 2\left(\frac{8}{21}\right)^x + \left(\frac{4}{11}\right)^x - 3\left(\frac{16}{45}\right)^x + \left(\frac{8}{25}\right)^x + \left(\frac{4}{13}\right)^x \ldots$$

In comparison to the expansion of $1/(\varsigma(x)-1)$, this one has irregularly changed signs, one missing summand, namely 2(4/9)^x and one missing coefficient − in front of (8/27)^x.

Let $\beth(s) = \frac{1}{2^s} + \sum_{\substack{n \in A \\ n > 2}} n^{-s} - \sum_{\substack{m \in B \\ m > 2}} m^{-s}$

We'll have, for large $s$:

$$\frac{1}{\beth(s)} = 2^s + \left(\frac{4}{3}\right)^s - 1 + \left(\frac{8}{9}\right)^s + \left(\frac{4}{5}\right)^s - 3\left(\frac{2}{3}\right)^s + \left(\frac{16}{27}\right)^s - \left(\frac{4}{7}\right)^s + 2\left(\frac{8}{15}\right)^s + 2\left(\frac{1}{2}\right)^s \ldots$$

Another example is provided by $1/(1-1/\varsigma(s)) = \varsigma(s)/(\varsigma(s)-1)$

$$\frac{1}{\zeta(s)} = \sum_{n=1}^{\infty} \frac{\mu(n)}{n^s}$$

(where $\mu$ is the aperiodic Möbius function. The expansion of $\varsigma(s)/(\varsigma(s) - 1)$ is exactly (♚) but without the constant term −1).

But let us turn back to Dirichlet L-series. By means of Dirichlet character

$$\chi_{-3}(n) = \begin{cases} 1 & \Longleftrightarrow \quad n \equiv 1 \quad (\text{mod } 3) \\ 0 & \Longleftrightarrow \quad n \equiv 0 \quad (\text{mod } 3) \\ -1 & \Longleftrightarrow \quad n \equiv 2 \quad (\text{mod } 3) \end{cases}$$

define the function: $\quad L(s, \chi_{-3}) = \sum_{n=1}^{\infty} \frac{\chi_{-3}(n)}{n^s} = 1 - \frac{1}{2^s} + \frac{1}{4^s} - \frac{1}{5^s} + \ldots$

Then, for large $s$

$$\frac{1}{1 - L(s, \chi_{-3})} = 2^s + 1 - \left(\frac{4}{5}\right)^s + \left(\frac{4}{7}\right)^s - \left(\frac{2}{5}\right)^s - \left(\frac{4}{11}\right)^s + \left(\frac{8}{25}\right)^s +$$

$$+ \left(\frac{4}{13}\right)^s + \left(\frac{2}{7}\right)^s - \left(\frac{4}{17}\right)^s - 2\left(\frac{8}{35}\right)^s + \left(\frac{4}{19}\right)^s - \left(\frac{2}{11}\right)^s - \left(\frac{4}{23}\right)^s$$

$$+ \left(\frac{8}{49}\right)^s + 2\left(\frac{4}{25}\right)^s + \ldots \tag{♝}$$

It seems that in the RHS of (♝) there are no denominators $\equiv 0$ (mod 3)… Those precisely which display this feature in (♚), are skipped in (♝).
One finds for

$$L(s, \chi_{-5}) = \frac{1}{1^s} - \frac{1}{2^s} - \frac{1}{3^s} + \frac{1}{4^s} + \frac{1}{6^s} - \frac{1}{7^s} - \frac{1}{8^s} + \frac{1}{9^s} + \frac{1}{11^s} - \ldots$$

the asymptotic expansion

$$\frac{1}{1 - L(s, \chi_{-5})} = 2^s - \left(\frac{4}{3}\right)^s + 1 + \left(\frac{8}{9}\right)^s - \left(\frac{2}{3}\right)^s - \left(\frac{16}{27}\right)^s - \left(\frac{4}{7}\right)^s + 2\left(\frac{4}{9}\right)^s +$$

$$+ \left(\frac{32}{81}\right)^s + 2\left(\frac{8}{21}\right)^s + \left(\frac{4}{11}\right)^s - \left(\frac{4}{13}\right)^s - 3\left(\frac{8}{27}\right)^s - \left(\frac{2}{7}\right)^s - \left(\frac{64}{243}\right)^s - 3\left(\frac{16}{63}\right)^s - 3\left(\frac{16}{63}\right)^s$$

$$-2\left(\frac{8}{33}\right)^s - \left(\frac{4}{17}\right)^s + \left(\frac{2}{9}\right)^s + \left(\frac{4}{19}\right)^s + 2\left(\frac{8}{39}\right)^s + 4\left(\frac{16}{81}\right)^s + \ldots \tag{☯}$$

Here again, there are no denominators divisible by 5. It seems that coefficients larger than 1 in absolute value show up often either when the denominators have at least two distinct factors or when the numerator and denominator are the same power of two different primes.

We have, for the two Dirichlet characters modulo 6, the beginnings of the two following expansions.

Let $$L(s,\chi_6) = 1 + \frac{1}{5^s} + \frac{1}{7^s} + \frac{1}{11^s} + \frac{1}{13^s} + \cdots$$

Then, for large $s$

$$\frac{1}{L(s,\chi_6) - 1} = 5^s - \left(\frac{25}{7}\right)^s + \left(\frac{125}{49}\right)^s - \left(\frac{25}{11}\right)^s - \left(\frac{25}{13}\right)^s - \left(\frac{625}{343}\right)^s +$$

$$+2\left(\frac{125}{77}\right)^s - \left(\frac{25}{17}\right)^s + 2\left(\frac{125}{91}\right)^s - \left(\frac{25}{19}\right)^s + \left(\frac{3125}{2401}\right)^s - 3\left(\frac{625}{539}\right)^s -$$

$$-\left(\frac{25}{23}\right)^s + 2\left(\frac{125}{119}\right)^s + \left(\frac{125}{121}\right)^s - 1 - 3\left(\frac{625}{637}\right)^s \cdots \qquad (\text{♟})$$

Let $$L(s,\chi_{-6}) = 1 - \frac{1}{5^s} + \frac{1}{7^s} - \frac{1}{11^s} + \frac{1}{13^s} - \cdots$$

Then, for large $s$

$$\frac{1}{1 - L(s,\chi_{-6})} = 5^s + \left(\frac{25}{7}\right)^s + \left(\frac{125}{49}\right)^s - \left(\frac{25}{11}\right)^s + \left(\frac{25}{13}\right)^s + \left(\frac{625}{343}\right)^s -$$

$$-2\left(\frac{125}{77}\right)^s - \left(\frac{25}{17}\right)^s + 2\left(\frac{125}{91}\right)^s + \left(\frac{25}{19}\right)^s + \left(\frac{3125}{2401}\right)^s - 3\left(\frac{625}{539}\right)^s -$$

$$-\left(\frac{25}{23}\right)^s - 2\left(\frac{125}{119}\right)^s + \left(\frac{125}{121}\right)^s + 1 + 3\left(\frac{625}{637}\right)^s - \cdots \qquad (\text{♙})$$

As one can see, (♟) and (♙) are similar up to some signs. The rationals the expansions (♟) and (♙) are made of are products of the small primes (with small exponents) which appear in the Dirichlet $L$-series $L(s,\chi_6)$ and $L(s,\chi_{-6})$: 7, 11, 13, 17, 19, etc.

Due to the structure of $\chi$, the numerators are always powers of 5; in turn, even as a factor, 3 does not appear in the denominators, where one can find the next primes: 7, 11, 13, 17, 19, etc. either as a factor or as an integer power of it (e.g. 343 = $7^3$ and 2401 = $7^4$). Examples of factorizations: 539=$7^2$*11, 637=$7^2$*13. The way new bigger primes appear in the denominators and then generate the next ones might be an interesting process. For example, the summands of rank 6, 7, 8 and 9 are products of the second summand – (25/7)^s – by, respectively, the summands of rank 3, 4, 5 divided by 5. Etc. etc. We'll end this chapter with the following remark: in (♚), which is about Riemann zeta function, the integer coefficients greater than 1 appear only near fractions that either have a composite denominator with at least two distinct prime factors or near a fraction whose numerator and denominator are exactly the same power (greater than 1) of 2 and, respectively, 3. In the known terms of (♚), if one neglects the constant term, the (odd) primes appear for the first time in the denominators exactly with their own density multiplied by 3/2. As will be shown in the next chapter, this frequency slows down after 13.

## 5.1. Comparing two characters modulo 10

For the two real Dirichlet characters modulo 10, the rationals related to $1/(L_{10}-1)$ are easily computable (up to a certain point).

For χ(1)=1, χ(3)=1, χ(7)=1, χ(9)=1 (otherwise 0), we have

$$\frac{1}{L(s,\chi_{10}) - 1} = 3^s - \left(\frac{9}{7}\right)^s - 1 - \left(\frac{9}{11}\right)^s - \left(\frac{9}{13}\right)^s + \left(\frac{27}{49}\right)^s - \left(\frac{9}{17}\right)^s$$

$$-\left(\frac{9}{19}\right)^s + \left(\frac{3}{7}\right)^s - \left(\frac{9}{23}\right)^s + 2\left(\frac{27}{77}\right)^s - \left(\frac{9}{29}\right)^s + 2\left(\frac{27}{91}\right)^s - \left(\frac{9}{31}\right)^s$$

$$+\left(\frac{3}{11}\right)^s - \left(\frac{9}{37}\right)^s - \left(\frac{81}{343}\right)^s + \left(\frac{3}{13}\right)^s + 2\left(\frac{27}{119}\right)^s + \left(\frac{27}{121}\right)^s - \left(\frac{9}{41}\right)^s$$

$$-\left(\frac{9}{43}\right)^s + 2\left(\frac{27}{133}\right)^s - \left(\frac{9}{47}\right)^s + 2\left(\frac{27}{143}\right)^s - 2\left(\frac{9}{49}\right)^s + \left(\frac{3}{17}\right)^s - \left(\frac{9}{53}\right)^s$$

$$+2\left(\frac{27}{161}\right)^s + \left(\frac{27}{169}\right)^s + \left(\frac{3}{19}\right)^s - \left(\frac{9}{59}\right)^s - 3\left(\frac{81}{539}\right)^s - \left(\frac{9}{61}\right)^s + 2\left(\frac{27}{187}\right)^s$$

$$-\left(\frac{9}{67}\right)^s + 2\left(\frac{27}{203}\right)^s + \left(\frac{3}{23}\right)^s + 2\left(\frac{27}{209}\right)^s - 3\left(\frac{81}{637}\right)^s - \left(\frac{9}{71}\right)^s + 2\left(\frac{27}{217}\right)^s$$

$$-\left(\frac{9}{73}\right)^s + 2\left(\frac{27}{221}\right)^s - 3\left(\frac{9}{77}\right)^s - \left(\frac{9}{79}\right)^s + \dots \quad (\Psi)$$

For the other real character modulo 10 (with χ(3) = χ(7) = –1 and χ(1) = χ(9) = 1), we have

$$\frac{1}{1 - L(s, \chi_{-10})} = 3^s - \left(\frac{9}{7}\right)^s + 1 + \left(\frac{9}{11}\right)^s - \left(\frac{9}{13}\right)^s + \left(\frac{27}{49}\right)^s - \left(\frac{9}{17}\right)^s$$

$$+ \left(\frac{9}{19}\right)^s - \left(\frac{3}{7}\right)^s - \left(\frac{9}{23}\right)^s - 2\left(\frac{27}{77}\right)^s + \left(\frac{9}{29}\right)^s + 2\left(\frac{27}{91}\right)^s + \left(\frac{9}{31}\right)^s$$

$$+ \left(\frac{3}{11}\right)^s - \left(\frac{9}{37}\right)^s - \left(\frac{81}{343}\right)^s - \left(\frac{3}{13}\right)^s + 2\left(\frac{27}{119}\right)^s + \left(\frac{27}{121}\right)^s + \left(\frac{9}{41}\right)^s$$

$$- \left(\frac{9}{43}\right)^s - 2\left(\frac{27}{133}\right)^s - \left(\frac{9}{47}\right)^s - 2\left(\frac{27}{143}\right)^s + 2\left(\frac{9}{49}\right)^s - \left(\frac{3}{17}\right)^s - \left(\frac{9}{53}\right)^s$$

$$+ 2\left(\frac{27}{161}\right)^s + \left(\frac{27}{169}\right)^s + \left(\frac{3}{19}\right)^s + \left(\frac{9}{59}\right)^s + 3\left(\frac{81}{539}\right)^s + \left(\frac{9}{61}\right)^s - 2\left(\frac{27}{187}\right)^s$$

$$- \left(\frac{9}{67}\right)^s - 2\left(\frac{27}{203}\right)^s - \left(\frac{3}{23}\right)^s + 2\left(\frac{27}{209}\right)^s - 3\left(\frac{81}{637}\right)^s \ldots \qquad (☬)$$

Time and again, we see integer powers of 3 in the numerators.
One can notice the sequence of primes in the denominators of fractions whose numerator is 9, and another sequence of ordered primes in fractions whose numerator is 3.

We never encounter numbers divisible by 5 in the denominator: one should remember that χ(5) = 0.
All composite numbers in the denominators which are not integer powers of primes bring about a coefficient equal to the numbers of its factors (taken with multiplicity).
An integer power of a prime does not give place to coefficients unless it coincides with the exponent in the numerator: e.g. 2(9/49), –2(3/23), 2(3/29)
The changes of signs obey a not at all obvious rule. In turn, the presence of coefficients greater than 1 (in absolute value) might be decrypted adding to the aforementioned rules the following one: if the number of distinct factors in the denominator coincides with the power of 3 in the numerator, then a coefficient greater than one shows up in front of the rational, e.g. –2(9/77).
Another comparison concerns the Riemann zeta function and the Dirichlet eta function. They belong to different moduli, the latter being the alternate version of the other.

## 5.2. Comparing zeta to eta

For the Riemann zeta function, the sequence of rationals is (with the coefficients written in parenthesis):
2, –4/3, –1, 8/9, –4/5, 2/3, –16/27, –4/7, (2) 8/15, (2) –4/9, 2/5, 32/81, (2) 8/21, –4/11, (3) –16/45, 8/25, –4/13, (3) 8/27, 2/7, (3) –4/15, –64/243, (3) –16/63, (2) 8/33, (4) 32/135, –4/17, (2) 8/35, 2/9, (3) –16/75, –4/19, (2) 8/39, (4) –16/81, (3) –4/21, 2/11, (7) 8/45,… (These may be found in (♚), here with several additional terms.)

For the Dirichlet eta function, the asymptotic expansion reads:

$$\frac{1}{1-\eta(x)} = 2^x + \left(\frac{4}{3}\right)^x - 1 + \left(\frac{8}{9}\right)^x + \left(\frac{4}{5}\right)^x - 3\left(\frac{2}{3}\right)^x + \left(\frac{16}{27}\right)^x + \left(\frac{4}{7}\right)^x$$

$$+2\left(\frac{8}{15}\right)^x - 4\left(\frac{4}{9}\right)^x - 3\left(\frac{2}{5}\right)^x + \left(\frac{32}{81}\right)^x + 2\left(\frac{8}{21}\right)^x + \left(\frac{4}{11}\right)^x + 3\left(\frac{16}{45}\right)^x + 2\left(\frac{1}{3}\right)^x$$

$$+\left(\frac{8}{25}\right)^x + \left(\frac{4}{13}\right)^x - 5\left(\frac{8}{27}\right)^x - 3\left(\frac{2}{7}\right)^x - 9\left(\frac{4}{15}\right)^x + \left(\frac{64}{243}\right)^x + 3\left(\frac{16}{63}\right)^x$$

$$+2\left(\frac{8}{33}\right)^x + 4\left(\frac{32}{135}\right)^x + \left(\frac{4}{17}\right)^x + 2\left(\frac{8}{35}\right)^x + 5\left(\frac{2}{9}\right)^x + 3\left(\frac{16}{75}\right)^x + \left(\frac{4}{19}\right)^x$$

$$+2\left(\frac{8}{39}\right)^x + 2\left(\frac{1}{5}\right)^x - 6\left(\frac{16}{81}\right)^x - 9\left(\frac{4}{21}\right)^x - 3\left(\frac{2}{11}\right)^x - 17\left(\frac{8}{45}\right)^x \quad (♚)$$

In spite of all the similitudes, in these expansions there are a lot of new items in comparison to the expansion of $1/(\varsigma(s) - 1)$. The summand $2(1/3)^x$ does not appear in the zeta expansion at all. Nor does $2(1/5)^x$. So far, we didn't find 1 in numerators of any other expansion of a Dirichlet *L*-function.

Here the coefficients are not at all the same as in (♚).
Interesting seems to be, in the eta Dirichlet asymptotic expansion, the presence of summands of the form

$$-(k+2)\left(\frac{2^k}{3^k}\right)^s \quad \text{or of the forms} \quad -3\left(\frac{2}{p}\right)^s \quad \text{and} \quad -(2^{k+1}+1)\left(\frac{2^k}{p_1\ldots p_k}\right)^s$$

If the denominator of a fraction is a composite number with at least two distinct factors, then the coefficient of that fraction is greater than 1. The coefficients appear to be particularly large when the denominator is not a power of a prime and when, at the same time, its number of factors (taken with multiplicity) coincides with the power of 2 in the numerator. Yet it is not clear wether and in which way it depends on the number of factors (distinct or not) in the denominator and/or on the power in the numerator. If the sign of a summand[1] is +, then coefficients coincide with the number of factors of the denominator, provided there exist at least two distinct factors.

If the number of factors (distinct or not) in the denominator matches the power or 2 in the numerator, then a coefficient > 1 (in absolute value) appears.

It is not clear whether any odd number divides at least one denominator in the presumably infinite (possibly displaying farther some fractional coefficients!) with sequence of fractions attached to a Dirichlet *L*-function.

It is not clear if the sets of denominators are or not closed under multiplication. Nor if the following weaker statement holds or not: for any denominators *n* and *m* in the set $D_{L(s,\chi)}$ of the denominators in the expansion of $1/(L(s,\chi)-1)$ there is a *k* in $D_{L(s,\chi)}$ so that *nm* divides *k*.

Statistics and dynamics of the sets of these rationals might reveal interesting.

An informal question: suppose you wake up tomorrow morning in a world where the correctness and the proof of GRH is not anymore a matter of discussion. How slightly an 'artificially' constructed or altered function may deviate from having all zeros on its 'critical line'? And which would be the criteria of 'slightness'?

## 6. Rationals and Hurwitz zeta function

Hurwitz zeta function is defined as

$$\zeta(s,q) = \sum_{n=0}^{\infty} \frac{1}{(q+n)^s}$$

We'll just put up two special cases in order to show the possibility of the study of rationals in this broader context:

$$\frac{1}{\zeta\left(x,\frac{3}{2}\right)} = \left(\frac{3}{2}\right)^x - \left(\frac{9}{10}\right)^x - \left(\frac{9}{14}\right)^x + \cdots \quad \text{(for large } x\text{)}$$

$$\frac{1}{\zeta\left(\frac{\pi}{2},x\right)} = \left(\frac{\pi}{2}\right)^x - \left(\frac{\pi^2}{2\pi+4}\right)^x + \cdots \quad \text{(also for large } x \text{ values)}$$

---

[1] this happens exactly in the cases when the Liouville λ function takes different values in the numerator and in the denominator (for a more general discussion see below)

One can of course define Hurwitz analogues of Dirichlet functions, although these would lose some of their essential properties (as does the Riemann zeta function itself). For example the Hurwitz lambda function would be:

$$\lambda(s, q) = \sum_{n=0}^{\infty} \frac{1}{(2n+q)^s}$$

One has for large values of the first argument

$$\frac{1}{\lambda\left(x, \frac{3}{2}\right)} = \left(\frac{3}{2}\right)^x - \left(\frac{9}{22}\right)^x - \ldots$$

and, respectively,

$$\frac{1}{\lambda\left(x, \frac{\pi}{2}\right)} = \left(\frac{\pi}{2}\right)^x - \left(\frac{\pi^2}{2\pi + 16}\right)^x - \ldots$$

## 7. On a connection with the Liouville lambda function

It would be nice if a connection with the Liouville function could be established at least for the asymptotic expansions of Dirichlet *L*-functions whose asymptotic part of the graph for real arguments lies above the line *y* = 1

The Liouville function, denoted by λ(*n*), is defined, for positive integers, as

$$\lambda(n) = (-1)^{\Omega(n)}$$

where Ω(*n*) is the number of prime factors of *n*, counted with multiplicity (sequence A008836 in OEIS)

For rationals, the Liouville function of the first kind, denoted by λ$_V$(*p/q*), might be defined as follows:

$$\lambda_V\left(\frac{p}{q}\right) = \begin{cases} \lambda(q) & \text{if } |\Omega(p) - \Omega(q)| = 1 \\ -\lambda(q) & \text{if } \Omega(p) = \Omega(q) \\ 1 & \text{if } p = 1 \text{ and } q \neq 1 \text{ (regardless of their } \Omega \text{ values)} \end{cases} \quad (\ast)$$

In fact, at least in (♚), (♙), (Ψ), (♖) (i.e. in the asymptotic expansions for Riemann zeta function, Dirichlet lambda, Dirichlet $L(s,\chi_6)$ and Dirichlet $L(s,\chi_{10})$ functions), the signs of the summands might (conjecturally) be rewritten exactly in terms of $\lambda_V(p/q)$. One can find additional heuristic confirmation in the beginning of the following asymptotic expansion:

$$\frac{1}{L(s,\chi_{+3})-1} = 2^s - 1 - \left(\frac{4}{5}\right)^s - \left(\frac{4}{7}\right)^s + \left(\frac{2}{5}\right)^s - \left(\frac{4}{11}\right)^s + \left(\frac{8}{25}\right)^s -$$

$$\left(\frac{4}{13}\right)^s + \left(\frac{2}{7}\right)^s - \left(\frac{4}{17}\right)^s + 2\left(\frac{8}{35}\right)^s - \left(\frac{4}{19}\right)^s + \left(\frac{2}{11}\right)^s - \left(\frac{4}{23}\right)^s +$$

$$+\left(\frac{8}{49}\right)^s - 2\left(\frac{4}{25}\right)^s + \left(\frac{2}{13}\right)^s + 2\left(\frac{8}{55}\right)^s + \ldots$$

The signs in the expansions of Dirichlet $L$-functions (♛), (♝), (●), (♙), (♗), whose real asymptotic part of the graphic lies under the line $y = 1$ (with the important exception of the Dirichlet eta function) seem to obey a completely different law: the signs seem to be partly determined by the residues of the denominators taken with respect to the Dirichlet character modulus.

More terms are needed to interpret heuristically the terms of the form $k*(2^k/3^k)^s$ in (●). The signs of the summands in the expansion (♚) related to the Dirichlet eta function deserves a separate discussion:

Here we define, for rationals, the Liouville function of the second kind as follows:

$$\lambda_A\left(\frac{p}{q}\right) = \begin{cases} 1 & \text{if } |\Omega(p) - \Omega(q)| = 1 \\ -1 & \text{if } \Omega(p) = \Omega(q) \end{cases} \qquad (☺)$$

where $\Omega(p)$ denotes again the number of prime factors of $p$ counted with multiplicity. (☺) completely characterizes the signs of known terms of (♚).

### NOTE

Both (※) and (☺) suffice to characterize the signs of the summands we were able to compute with sufficient accuracy. Both these formulations are liable to be refined and/or completed when additional terms of the expansions of $1/(\varsigma(x)-1)$ and of $1/(1-\eta(1))$ will be properly calculated. They have the signs predicted by, respectively, (※) and (☺), but we have to be cautious: for example, the case $|\Omega(p) - \Omega(q)| > 1$ remains unclear: until now, we haven't encountered, in the studied expansions, any term of this kind whatsoever.

It is not yet clear whether the summands $+2(1/3)^s$ and $+2(1/5)^s$ are positive because $|\Omega(1) - \Omega(5)| = 1$ and $|\Omega(1) - \Omega(3)| = 1$ or just because, as in (♚), it should be supposed, as for the Riemann zeta function (see (※)), that $\lambda_A(1/q) = 1$ for all $q$. At least heuristically, these questions will rapidly find their answers, provided one might be interested in further computation and study of these expansions.

Liouville lambda function is related to Riemann zeta function by the well-known expansion formula:

$$\frac{\zeta(2s)}{\zeta(s)} = \sum_{n=1}^{\infty} \frac{\lambda(n)}{n^s} \qquad (♣) \qquad \text{(by the way, another aperiodic character)}$$

Interestingly, the expansion is made of exactly the same summands as the corresponding to Riemann zeta function expansion (♚) up to some signs.

$$\frac{1}{1 - \frac{\zeta(2s)}{\zeta(s)}} = \frac{\zeta(s)}{\zeta(s) - \zeta(2s)} = 2^s - \left(\frac{4}{3}\right)^s + 1 + \left(\frac{8}{9}\right)^s - \left(\frac{4}{5}\right)^s - \left(\frac{2}{3}\right)^s - \left(\frac{16}{27}\right)^x$$

$$- \left(\frac{4}{7}\right)^s + 2\left(\frac{8}{15}\right)^s + 2\left(\frac{4}{9}\right)^s - \left(\frac{2}{5}\right)^s + \left(\frac{32}{81}\right)^s + 2\left(\frac{8}{21}\right)^s - \left(\frac{4}{11}\right)^s$$

$$-3\left(\frac{16}{45}\right)^s + \left(\frac{8}{25}\right)^s - \left(\frac{4}{13}\right)^s - 3\left(\frac{8}{27}\right)^s - \left(\frac{2}{7}\right)^s + 3\left(\frac{4}{15}\right)^s - \left(\frac{64}{243}\right)^s$$

$$-3\left(\frac{16}{63}\right)^s + 2\left(\frac{8}{33}\right)^s + 4\left(\frac{32}{135}\right)^s - \left(\frac{4}{17}\right)^s + 2\left(\frac{8}{35}\right)^s + \left(\frac{2}{9}\right)^s - 3\left(\frac{16}{75}\right)^s$$

$$-\left(\frac{4}{19}\right)^s + 2\left(\frac{8}{39}\right)^s + 4\left(\frac{16}{81}\right)^s + 3\left(\frac{4}{21}\right)^s - \left(\frac{2}{11}\right)^s - 7\left(\frac{8}{45}\right)^s \ldots \quad (❀)$$

The signs of this asymptotic expansion (❀) are directly related to the Liouville $\lambda$-function applied only to the denominators regardless of the number of prime factors (i.e. regardless of the exponent of 2) in the numerator.

Yet, (♣) holds for any $s > 1$, while the asymptotical (❀) holds only in some neighborhood of infinity.

This is a naive but not necessarily false way of thinking that primes 'behave at infinity' basically in the way they behave anywhere. Anyhow, since $s$ tends to infinity, tiny primes and prime factors of larger numbers (e.g. $243 = 3^5$, $135 = 5*3^3$ ) are mirrored in the 'neighborhoods of infinity' through a function closely related to the Riemann zeta function, namely $1/(\varsigma(s-1))$.

Andrei Vieru

andreivieru007@gmail.com